\documentclass[10pt]{article}
\usepackage{amsfonts}

\usepackage{amsmath}
\usepackage{graphicx}

\newtheorem{theorem}{Theorem}[section]

\newtheorem{lemma}[theorem]{Lemma}

\newtheorem{problem}[theorem]{Problem}

\newenvironment{proof}[1][Proof]{\textbf{#1.} }{\ \rule{0.5em}{0.5em}}

\def\si{{\sigma}}

\def\al{{\alpha}}

\def\Ga{{\Gamma}}

\def\si{{\sigma}}

\def \eref#1{\hbox{(\ref{#1})}}
\def\th{{\theta}}

\def \eref#1{\hbox{(\ref{#1})}}
\def\th{{\theta}}
\def\si{{\sigma}}
\def\al{{\alpha}}

\begin{document}

\title{Parameter estimation for  fractional Ornstein-Uhlenbeck processes}
\author{Yaozhong Hu\thanks{%
Y. Hu is supported by the National Science Foundation under DMS0504783 }%
\quad and\quad D. Nualart\thanks{%
D. Nualart is supported by the National Science Foundation under DMS0604207 }
\\
Department of Mathematics\thinspace ,\ University of Kansas\\
405 Snow Hall\thinspace ,\ Lawrence, Kansas 66045-2142\\
hu@math.ku.edu and nualart@math.ku.edu}
\date{}
\maketitle

\begin{abstract}
We   study a least squares estimator $\widehat {\theta}_T$ for the
Ornstein-Uhlenbeck process, $dX_t=\theta X_t dt+\sigma dB^H_t$,
 driven by fractional Brownian motion $B^H$ with Hurst
parameter $H\ge \frac12$. We prove the strong consistence 
of $\widehat {\theta}_T$ (the almost surely convergence 
of $\widehat {\theta}_T$ to the true parameter ${%
\theta}$). 
We also obtain the rate of this convergence when $1/2\le H<3/4$, 
applying a central limit theorem for multiple Wiener integrals.  This
least squares estimator can be used to study  other more simulation friendly 
estimators such as the estimator $\tilde \th_T$ defined by \eref{g11}. 
\end{abstract} 

\section{Introduction}

The Ornstein-Uhlenbeck process $X_t$ driven by a certain type of noise $Z_t$
is described by the Langevin equation
\begin{equation*}
  X_t=X_0- \theta  \int_0^tX_s ds+\sigma  Z_t.
\end{equation*}
If the parameter $\theta$ is unknown and if the process $(X_t, 0\le t\le T)$
can be observed continuously, then an important problem is to estimate the
parameter ${\theta }$ based on the (single path) observation $(X_t, 0\le t\le
T)$. When $Z_t$ is the standard Brownian motion, this problem has been
extensively studied (see for example \cite{Ku}, \cite{LS} and the references
therein). The most popular approaches are either the maximum likelihood
estimators or the least squares estimators, and in this case
they coincide. Other type of noise processes have also been studied. For
example, when $Z_t$ is an ${\alpha}$-stable process  maximum likelihood
estimators do  not exist and other approaches are proposed in \cite{HL} and \cite{HL1}.

In this paper we study the parameter estimation problem for the  Ornstein-Uhlenbeck process driven by fractional Brownian motion  with Hurst parameter $H $
\begin{equation}
 X_t=X_0- \theta \int_0^t X_sds+\sigma  B_t^H,   \label{e.1.1}
\end{equation}
where ${\theta }>0$ is an unknown parameter.
Although the Ornstein-Uhlenbeck  process is defined for all $H\in (0, 1)$, we assume
$H>1/2$ in this paper.
  In \cite{KL}, the  the maximum likelihood estimator $\bar {\theta}_T $ for the parameter $\theta$ is obtained and
   has the following expression
\begin{equation*}
\bar {\theta}_T=-\left\{ \int_0^T Q^2(s) dw_s^H\right\}^{-1} \int_0^T Q(s)
dZ_s\,,
\end{equation*}
where
\begin{eqnarray*}
k_H(t,s)&=&\kappa_H ^{-1} s^{\frac 12-H}(t-s)^{\frac 12-H}\,,\quad \kappa_H=2H
\Gamma\left(\frac 32-H\right)\Gamma\left(H+\frac 12\right)\,; \\
w_t^H&=&{\lambda}_H^{-1} t^{2-2H}\,;\quad {\lambda}_H=\frac{2H{\Gamma}(3-2H){%
\Gamma}\left(H+\frac 12\right)}{{\Gamma}\left(\frac 32-H\right)}\,; \\
Q(t)&=&\frac{d}{dw_t^H} \int_0^t k_H(t,s) X_s ds\,,\quad 0\le t\le T\,; \\
Z_t&=&\int_0^t k_H(t,s) dX_s\,.
\end{eqnarray*}
It is proved
that  $\lim_{T\rightarrow \infty} \bar {\theta}_T={\theta}$ almost surely.

In this paper we   propose   two different estimators for the parameter $\theta$ and we study their asymptotic behavior. First we introduce an estimator of the form
\begin{equation}  \label{eq1}
\widehat \theta_T=  \theta - {\sigma}
\frac{\int_0^T X_t   dB^H_t}{\int_0^T X_t^2 dt},
\end{equation}
where  $\int_0^T X_t   dB^H_t $ is a divergence-type integral (see \cite{BHOZ}, \cite{DHP}, \cite{Hu}, \cite{HO}  and the references therein), and we call it
the {\it  least squares estimator}.  This is motivated by the following heuristic argument.
The least square  estimator  aims to minimize
\begin{equation*}
\int_0^T |\dot X_t+\theta X_t|^2 dt\,,
\end{equation*}
and this leads to the solution
\begin{equation}  \label{e.1.3}
\widehat \theta_T=-\frac{\int_0^T X_t dX_t}{\int_0^T X_t^2 dt}.
\end{equation}
If $H=\frac 12$, then   the integral  $\int_0^T X_t dX_t$ is an It\^o stochastic integral which can be approximated by forward Riemann sums. However, for $H>\frac 12$  the numerical simulation of the estimator  $\widehat \theta_T$ seems extremely difficult. 
 For this reason, in this case we introduce  and study 
 a  second  estimator $\widetilde \theta_T$, defined in (\ref{g11}).

We prove the almost sure convergence of the  estimator $\widehat\theta_T$ to $\theta$, as $T$ tends to infinity, and  derive  the rate of convergence,   obtaining a central limit theorem in the case $H\in \left[\frac 12 ,\frac 34\right)$.  The proof of the central limit theorem is based on the characterization of the convergence in   law for multiple stochastic integrals using the techniques of Malliavin calculus,  established recently by Nualart and Ortiz-Latorre in
 \cite{NO}.     Finally, we derive  the rate of convergence of the estimator  $\widetilde \theta_T$ from the rate of convergence of  $\widehat \theta_T$.

\setcounter{equation}{0}

\section{Preliminaries}

In this section we first introduce some basic facts on the Malliavin calculus for the fractional Brownian motion and recall the main
result in \cite{NO} concerning the central limit theorem for multiple
stochastic integrals.

The fractional Brownian motion with Hurst
parameter $H\in (0,1)$,  $(B_{t}^{H},t\in  \mathbb{R})$   is a zero mean  Gaussian
process with covariance
\begin{equation}
\mathbb{E}(B_{t}^{H}B_{s}^{H})=R_{H}(s,t)=\frac{1}{2}\left(
|t|^{2H}+|s|^{2H}-|t-s|^{2H}\right) \,.  \label{cov}
\end{equation}%
We assume that $B^H$ is defined on a complete probability space $(\Omega,\mathcal{A},P)$ such that $\mathcal{A}$ is generated by $B^H$.
Fix a time interval $[0,T]$. Denote by $\mathcal{E}$ the set of real valued
step functions on $[0,T]$ and let $\mathcal{H}$ be the Hilbert space defined
as the closure of $\mathcal{E}$ with respect to the scalar product
\[
\langle
\mathbf{1}_{[0,t]},\mathbf{1}_{[0,s]}\rangle _{\mathcal{H}}=R_{H}(t,s),
\]
where $R_{H}$ is the covariance function of the fBm, given in (\ref{cov}).
The mapping $\mathbf{1}_{[0,t]}$$\longmapsto B_{t}^{H}$ can be extended to
a linear isometry between $\mathcal{H}$ and the Gaussian space $\mathcal{H}_{1}$
spanned by $B^{H}$. We denote this isometry by $\varphi \longmapsto
B^{H}(\varphi )$. For $H=\frac 12$ we have $\mathcal{H}= L^2([0,T])$, whereas
for $H>\frac 12$ we have $L^{\frac 1H} ([0,T])\subset \mathcal{H}$  
and for   $\varphi, \psi \in L^{\frac 1H} ([0,T])$ we have
\begin{equation} \label{eq9}
 \langle \varphi, \psi \rangle_{\mathcal{H}}= \alpha_H\int_0^T\int_0^T 
  \varphi_s\psi_t  |t-s|^{2H-2}dsdt,
\end{equation}
where $\alpha_H=H(2H-1)$.

Let $\mathcal{S}$ be the space of smooth and cylindrical random variables of
the form%
\begin{equation}
F=f(B^{H}(\varphi _{1}),\ldots ,B^{H}(\varphi _{n})),  \label{g1}
\end{equation}%
where $f\in C_{b}^{\infty }(\mathbb{R}^{n})$ ($f$ and all its partial
derivatives are bounded). For a random variable $F$ of the form \ (\ref{g1})
we define its Malliavin derivative as the $\mathcal{H}$-valued random
variable%
\begin{equation*}
DF=\sum_{i=1}^{n}\frac{\partial f}{\partial x_{i}}(B^{H}(\varphi
_{1}),\ldots ,B^{H}(\varphi _{n}))\varphi _{i}\text{. }
\end{equation*}%
 By iteration, one can
define the $m$th derivative $D^{m}F$, which is an element of $L^{2}(\Omega ;%
\mathcal{ H}^{\otimes m})$,
for every $m\geq 2$.
For $m\geq 1$, ${\mathbb{D}}^{m,2}$ denotes the closure of $%
\mathcal{S}$ with respect to the norm $\Vert \cdot \Vert _{m,2}$, defined by
the relation
\begin{equation*}
\Vert F\Vert _{m,2}^{2}\;=\;\mathbb{E}\left[ |F|^{2}\right] +\sum_{i=1}^{m}\mathbb{E}\left(
\Vert D^{i}F\Vert _{\mathcal{ H}^{\otimes i}}^{2}\right) .
\end{equation*}%
Let $\delta $ be the adjoint of the operator $D$, also called the
\textsl{divergence operator}. A random element $u\in L^{2}(\Omega ,\mathcal{H} ) $ belongs to the domain of $\delta $, denoted $\mathrm{Dom}(\delta)     $, if and
only if it verifies
\begin{equation*}
|E\langle DF,u\rangle _{\mathcal{ H}}|\leq c_{u}\,\Vert F\Vert _{L^{2}},
\end{equation*}%
for any $F\in \mathbb{D}^{1,2}$, where $c_{u}$ is a constant depending only
on $u$. If $u\in \mathrm{Dom}(\delta)     $, then the random variable $\delta (u)$
is defined by the duality relationship
\begin{equation}
E(F\delta (u))=E\langle DF,u\rangle _{\mathcal{ H}},  \label{ipp}
\end{equation}%
which holds for every $F\in {\mathbb{D}}^{1,2}$.    The divergence operator $\delta $ is also called the
Skorohod integral because in the case of the Brownian motion it coincides
with the anticipating stochastic integral introduced by Skorohod in \cite{Sk}.   We will make use of the notation $\delta(u)=\int_0^T u_t dB^H_t$.

For every $n\geq 1$, let $\mathcal{H}_{n}$ be the $n$th Wiener chaos
of $B,$
that is, the closed linear subspace of $L^{2}\left( \Omega ,\mathcal{A}%
,P\right) $ generated by the random variables $\{H_{n}\left( B^H\left(
h\right) \right) ,h\in H, \| h\| _{ \mathcal{ H}}=1\}$, where $H_{n}$ is
the $n$th Hermite polynomial. The mapping $I_{n}(h^{\otimes
n})=n!H_{n}\left( B^H\left( h\right) \right) $ provides a linear isometry
between the symmetric tensor product $\mathcal{ H}^{\odot n}$ and $\mathcal{H}_{n}$.
For $H=\frac{1}{2}$, $I_{n}$ coincides with the multiple It\^o
stochastic
integral.  On the other hand, $I_n(h^{\otimes n})$ coincides with the iterated divergence
$\delta^n(h^{\otimes n})$.

We will make use of the following central limit theorem for multiple  stochastic integrals
(see \cite{NO}).

\begin{theorem}
\label{t.2.2} Let $\{F_{n}\,,n\geq 1\}$ be a sequence of random variables in
the $p$-th Wiener chaos, $p\geq 2$, such that $\lim_{n\rightarrow \infty }%
\mathbb{E}(F_{n}^{2})=\sigma ^{2}$. Then the following conditions are
equivalent:
\begin{itemize}
\item[(i)] $F_n$ converges in law to $N(0,\sigma^2)$ as $n$ tends to
infinity.
\item[(ii)] $ \|DF_n\|^2_{\mathcal{H}} $ converges in $L^2$ to a constant   as  $n$ tends to infinity.
\end{itemize}
\end{theorem}

\noindent\textbf{Remark.}
In \cite{NO} it is proved that (i) is equivalent to the fact that $\Vert
DF_{n}\Vert _{\mathcal{H}}^{2}$ converges in $L^{2}$ to $p\sigma ^{2}$ as $n$ tends to
infinity.  If we assume (ii),    the limit of  $\|DF_n\|^2_{\mathcal{H}}$     must be equal to  $p\sigma ^{2}$ because
\[
\mathbb{E} (\|DF_n\|^2_{\mathcal{H}}) = p \mathbb{E}(F_{n}^{2}).
\]

\setcounter{equation}{0}

\section{Asymptotic behavior of the least square estimator}

Consider   Equation (\ref{e.1.1}) driven by a  fractional Brownian motion $B^H$ with Hurst parameter $H\ge \frac 12$. Suppose that $X_0=0$ and $\theta >0$.
The solution is given by
\begin{equation}
X_{t}=\sigma \int_{0}^{t}e^{-\theta (t-s)}dB^H_{s},
\label{a2}
\end{equation}%
where the stochastic integral is an It\^o integral if $H=\frac 12$ and a path-wise Riemann-Stieltjes integral if $H>\frac 12$.
Let  $\widehat \theta_T$ be the least squares estimator  defined in (\ref{eq1}).
The next lemma provides a useful alternative expression for    $\widehat \theta_T$.

\begin{lemma}
Suppose that $H>\frac 12$. Then
\begin{equation}  \label{e.3.7}
 \widehat \theta_T=  - \frac{X_T^2}{2\int_0^T X_t^2 dt}
+\sigma^2 \frac {\alpha_H \int_0^T \int_0^t \xi ^{2H-2} e^{{
-\theta} \xi } d\xi  dt}{\int_0^T X_t^2 dt}.
\end{equation}
\end{lemma}

\noindent
\textit{Proof} \quad
 Using the relation between the divergence integral and the path-wise Riemann-Stieltjes integral (see Theorem 3.12 and Equation (3.6) of \cite{DHP})  we can write
\begin{eqnarray*}
\int_0^T X_t  \circ d B_t^H
&=&\int_0^T X_t d B_t^H+  \alpha_H \int_0^T\int_0^t  D_sX_t (t-s) ^{2H-2}   ds dt  \\
&=&\int_0^T X_t d B_t^H+ \sigma  \alpha_H \int_0^T\int_0^t   e^{-\theta(t-s)} (t-s) ^{2H-2}   ds dt \\
&=& \int_0^T X_t d B_t^H+{\sigma}\alpha_H \int_0^T \int_0^t \xi ^{2H-2} e^{{
-\theta} \xi } d\xi  dt \,.
\end{eqnarray*}
As a consequence, we obtain
\begin{equation}  \label{e1}
\widehat \theta_T= \theta-  \sigma \frac{\int_0^T X_t  \circ d B_t^H}{\int_0^T X_t^2 dt}
+\sigma^2 \frac {\alpha_H \int_0^T \int_0^t \xi ^{2H-2} e^{{
-\theta} \xi } d\xi  dt}{\int_0^T X_t^2 dt}.
\end{equation}
On the other hand,
\begin{equation}  \label{e2}
\sigma \int_0^T X_t  \circ d B_t^H=   \int_0^T X_t  \circ d X_t + \theta\int_0^T X_t^2 dt
=\frac 12 X_T^2+ \theta\int_0^T X_t^2 dt.
\end{equation}
Substituting (\ref{e2}) into (\ref{e1}) yields (\ref{e.3.7}). 
\hfill $\Box$

 \medskip
The next theorem establishes the strong consistency of this estimator.

\begin{theorem}  \label{T1}
If $  H\ge \frac 12 $, then
\begin{equation}  \label{e3}
\widehat{\theta}_{T}\rightarrow {\theta }
\end{equation}%
almost surely, as $T$ tends to infinity.
\end{theorem}

In order to prove this theorem we make use of the following technical result.

\begin{lemma}  \label{L1}
Assume  $  H \ge \frac 12$. Then,
\begin{equation}
\frac{1}{T}\int_{0}^{T}X_{t}^{2}dt\rightarrow \sigma^2 \theta^{-2H}H\Gamma(2H),
\label{e.2.1}
\end{equation}%
almost surely and in $L^2$, as $T$ tends to infinity.
\end{lemma}

 \medskip
 \noindent
\textit{Proof} \quad
For every $t\ge 0$ define
\begin{equation}  \label{f6}
Y_t=\sigma \int_{-\infty }^{t}e^{-\theta
(t-s)}dB^H_{s}=X_{t}+e^{-\theta t}\xi,
\end{equation}%
where $\xi= \sigma \int_{-\infty}^0 e^{\theta s} dB^H_s$. The stochastic process
$(Y_t,t\geq 0)$ is Gaussian, stationary and ergodic. For $H=\frac 12$ this is well-known and for $H>\frac 12$ this is proved in  \cite{CKM}. Then, the ergodic theorem implies that
\begin{equation*}
\frac{1}{T}\int_{0}^{T}Y_{t}^2dt\rightarrow  \mathbb{E} (Y_0^2),
\end{equation*}
as $T$ tends to infinity, almost surely and in $L^2$.
This implies that
\begin{equation*}
\frac{1}{T}\int_{0}^{T}  X _{t}^2dt\rightarrow  \mathbb{E} (Y_0^2),
\end{equation*}
as $T$ tends to infinity, almost surely and in $L^2$.
If $H=\frac 12$, we know that  $ \mathbb{E} (Y_0^2)=  \frac {\sigma^2} { 2\theta}$, which implies  (\ref{e.2.1}).  If $H>\frac 12$,  using  (\ref{eq9}) yields
\[
\mathbb{E} (X_0^2)=\alpha_H \sigma^2 \int_0^\infty \int_0^\infty e^{-\theta(s+u)}   |u-s|^{2H-2} duds,
\]
and (\ref{e.2.1}) follows from Lemma \ref{A1}.
\hfill $\Box$

 \medskip
 \noindent
\textit{Proof  of Theorem  \ref{T1}} \quad
In the case $H=\frac 12$, taking into account that the process $\left(
\int_{0}^{t}X_{s}dB_{s},t\geq 0\right) $ is a martingale with quadratic
variation $\int_{0}^{t}X_{s}^{2}ds$ it follows that
$\widehat{\theta }_{T}\rightarrow \theta  $
almost surely, as $T$ tends to infinity.

Now let $H>1/2$.   From Lemma \ref{L6} we deduce that almost surely
\begin{equation}  \label{f5}
 \lim_{T\rightarrow \infty}\frac{X_T^2}T = 0.
\end{equation}
It is easy to check that this convergence also holds in $L^2$.
Then we conclude the proof using Lemma  \ref{L1}, (\ref{f5}), and
\[
\lim_{T\rightarrow \infty} \frac 1T   \int_0^T \int_0^t \xi ^{2H-2} e^{{
-\theta} \xi } d\xi  dt=  \theta^{1-2H} \Gamma(2H-1).
\]
\hfill $\Box$

The next theorem provides
 the convergence in distribution to a Gaussian law of the fluctuations in the almost sure
convergence (\ref{e3}).

\begin{theorem}\label{t.3.3}
Suppose $H\in \left[ \frac 12, \frac 34 \right)$. Let   $(X_{t}, t\in [0,T])$ be given by (\ref{a2}%
), then
\begin{equation}
\sqrt{T}\left[ \widehat{\theta }_{T}-\theta \right] \overset{\mathcal{L}}{%
\rightarrow }N(0,   \theta  \sigma_H^2)\ ,  \label{e.2.3}
\end{equation}%
as $T$ tends to infinity, where
\begin{equation}  \label{si}
\sigma^2_H= (4H-1)\left( 1+ \frac{\Gamma(3-4H)\Gamma(4H-1)}{\Gamma(2-2H)\Gamma(2H)}\right).
\end{equation}
\end{theorem}

\medskip
\noindent
\textit{Proof} \quad  We have
\begin{equation}
\widehat{{\theta }}_{T}-{\theta }=-{\sigma }\frac{\int_{0}^{T}X_{t}dB_{t}^{H}}{%
\int_{0}^{T}X_{t}^{2}dt}=-\frac{\sigma ^{2}\int_{0}^{T}\left( \int_{0}^{t}e^{{%
-\theta }(t-s)}dB_{s}^{H}\right) dB_{t}^{H}}{\int_{0}^{T}X_{t}^{2}dt}= -\frac{\sqrt{T}F_{T}}{ \int_{0}^{T}X_{t}^{2}dt},  \label{e.3.2}
\end{equation}%
where $F_T$ is the double stochastic integral
\begin{equation} \label{eq4}
F_{T}=   \frac {\sigma^2} {2\sqrt{T}}  I_2\left( e^{-\theta |t-s|}\right).
\end{equation}%
From Lemma \ref{L1} we know that $\frac 1T \int_{0}^{T}X_{t}^{2}dt$ converges almost surely and in $L^2$, as $T$ tends to infinity to  $\sigma^2 \theta^{-2H} H \Gamma(2H)$. Then, it suffices to show that  $F_T$ converges in law as $T$ tends to infinity to a  centered normal distribution. In order to show this convergence we will apply Theorem
\ref{t.2.2} to a given sequence of random variables in the second chaos $F_{T_k}$, where $T_k \uparrow \infty$ as $k$ tends to infinity. To simplify we assume that $T=1,2,\dots$. The proof then follows from the following facts:
\begin{itemize}
\item[(i)]  $\mathbb{E} (F_T^2)$ converges to  $\theta^{1-4H}\sigma^4  \delta_H$, where
\[
\delta_H=  H^2 (4H-1) (\Gamma (2H)^2 + \frac{\Gamma(2H)\Gamma(3-4H) \Gamma(4H-1)}{\Gamma(2-2H)}),
\]
 as $T$ tends to infinity.
\item[(ii)]  $   \| DF_T\|_{\mathcal{H}}^2  $ converges in $L^2$ to a constant  as $T$ tends to infinity.

\end{itemize}
\noindent
\textit{Step 1: Proof of (i)} \quad
Suppose first that $H=\frac 12$. In this case, by the isometry of the It\^o integral we obtain
\begin{equation*}
\mathbb{E} \left( F_{T}^2\right) =\frac{\sigma^4}{T}%
\int_{0}^{T}\int_{0}^{t}e^{-2\theta (t-s)}dsdt\ =\frac{\sigma^4}{T}\left(  %
\frac{T}{2{\theta }}+\frac{e^{-2{\theta }T}-1}{4{\theta }^{2}}\ \right),
\end{equation*}%
which implies that
\begin{equation*}
\lim_{T\rightarrow \infty }\mathbb{E} (F_{T}^{2})=  \frac{\sigma^4} {2\theta}.
\end{equation*}%
This implies the desired result because $\delta_{\frac 12}=\frac 12$.

Now, let $H\in \left(\frac 12, \frac 34 \right)$.
In this case, by the isometry property of the double stochastic integral $I_2$, the variance of $F_T$ is given by
\begin{equation}  \label{eq2}
 \mathbb{E} \left( F_{T}^{2}\right) =  \frac {\sigma^4 \alpha_H^2} {2T} I_T,
\end{equation}
where
\begin{equation}  \label{e5}
I_T= \int_{[0,T]^{4}}e^{{%
-\theta }|s_{2}-u_{2}|- { \theta }%
|s_{1}-u_{1}|}|s_{2}-s_{1}|^{2H-2}|u_{2}-u_{1}|^{2H-2}duds.
\end{equation}
By Lemma \ref{A2} in the Appendix we  obtain that
\[
\lim_{T\rightarrow \infty }\mathbb{E} (F_{T}^{2}) =\theta^{1-4H} \sigma^4 \delta_H.
\]

 \smallskip
\noindent
\textit{Step 2: Proof of (ii)} \quad
For $s\leq T$ we have
\begin{equation*}
D_{s}F_{T}=\frac{\sigma X_{s}}{\sqrt{T}}+\frac{\sigma^2 }{\sqrt{T}}%
\int_{s}^{T}e^{-\theta (t-s)}dB^H  _{t}.
\end{equation*}%
Suppose first that $H=\frac 12$. In this case,
\begin{eqnarray*}
\Vert DF_{T}\Vert _{\mathcal{H}}^{2} &=&\frac{\sigma^2 }{T}\int_{0}^{T}\left(
X_{s}^{2}+2\sigma X_{s}\int_{s}^{T}e^{-\theta (t-s)}dB_{t}+ \sigma^2\left(
\int_{s}^{T}e^{-\theta (t-s)}dB_{t}\right) ^{2}\right) ds \\
&=&A_{T}^{(1)}+A_{T}^{(2)}+A_{T}^{(3)}.
\end{eqnarray*}%
We already know from (\ref{e.2.1}) that $A_{T}^{(1)}$ converges in $L^{2}$ to $%
  \frac{\sigma^4} {2\theta}$ as $T$ tends to infinity. The third term can be
written as
\begin{equation*}
A_{T}^{(3)}=\frac{\sigma^4 }{T}\int_{0}^{T}\left( \int_{s}^{T}e^{-\theta
(t-s)}dB_{t}\right) ^{2}ds=\frac{\sigma^4 }{T}\int_{0}^{T}\left(
\int_{0}^{u}e^{-\theta (u-x)}dB_{x}\right) ^{2}du,
\end{equation*}%
so it also converges in $L^{2}$ to $  \frac{\sigma^4} {2\theta}$ a $T$ tends to
infinity. Finally we can show that
\begin{equation*}
\lim_{T\rightarrow \infty }\mathbb{E}\ ((A_{T}^{(2)})^{2})=0\,.
\end{equation*}%
In fact, we have
\begin{eqnarray*}
\mathbb{E}\ ((A_{T}^{(2)})^{2}) &=&\frac{8\sigma ^{6}}{T^{2}}\int_{\{s<u\leq
T\}}\mathbb{E}\ \left( X_{s}X_{u}\left( \int_{s}^{T}e^{-\theta
(t-s)}dB_{t}\right) \left( \int_{u}^{T}e^{-\theta (t-u)}dB_{t}\right) \right)
dsdu \\
&=&\frac{8\sigma ^{8}}{T^{2}}\int_{\{s<u\leq T\}}\left(
\int_{0}^{s}e^{-\theta (s+u-2r)}dr\right) \left( \int_{u}^{T}e^{-\theta
(2t-s-u)}dt\right) dsdu \\
&=&\frac{8\sigma ^{8}}{4\theta ^{2}T^{2}}\int_{\{s<u\leq T\}}( e^{2\theta
s}-1)( e^{-2\theta u}-e^{-2\theta T})dsdu,
\end{eqnarray*}%
which clearly converges to zero as $T$ tends to infinity. Therefore, $\Vert
DF_{T}\Vert _{\mathcal{H}}^{2}$ converges to  $  \frac{\sigma^4} { \theta}$ in $L^{2}$.

Suppose now that $H>\frac 12$.
 From (\ref{eq9}) we have
\begin{eqnarray*}
\Vert DF_{T}\Vert _{\mathcal{H}}^{2} &=&\frac{\alpha_{H}\sigma^2}{T}\int_{0}^{T}\int_{0}^{T}\left(
X_{s}+\si \int_{s}^{T}e^{-\theta (t-s)}dB_{t}^{H}\right)  \\
&&\times \left( X_{u}+\si \int_{u}^{T}e^{-\theta (t-u)}dB_{t}^{H}\right).
|u-s|^{2H-2}duds
\end{eqnarray*}%
We have to prove that  $\| DF_T \|^2_{\mathcal{H}}$ converges to a constant in $L^2$ as $T$ tends to infinity.  In fact,
\begin{eqnarray*}
\Vert DF_{T}\Vert _{\mathcal{H}}^{2} &=&\frac{ \alpha_H \sigma^2}{T}\int_{0}^{T}\int_{0}^{T}\Bigg(%
X_{s}X_{u}+2\si X_{u}\int_{s}^{T}e^{-\theta (t-s)}dB_{t}^{H} \\
&&+\si^2 \int_{s}^{T}e^{-\theta (t-s)}dB_{t}^{H}\int_{u}^{T}e^{-\theta
(t-u)}dB_{t}^{H}\Bigg)|u-s|^{2H-2}duds \\
&=&\frac{ \alpha_H \sigma^2}{T}( C_T^{(1)} + C_T^{(2)}+C_T^{(3)}).
\end{eqnarray*}%
For the term  $C_T^{(1)} $, since $X_t$ is Gaussian we can write
\begin{eqnarray*}
&& \mathbb{E}\left( |C_T^{(1)} - \mathbb{E}(C_T^{(1)}) |^2 \right)
=2  \int_{[0,T]^{4}}\mathbb{\ E}\ (X_{s}X_{t})\mathbb{\ E}\
(X_{u}X_{v})  \\
&&\quad \times |u-s|^{2H-2}|v-t|^{2H-2}dudvdsdt.
\end{eqnarray*}
By Lemma \ref{lem1}
\begin{eqnarray*}
 &&\frac{1}{T}\int_{[0,T]^{3}}\mathbb{\ E}\ (X_{T}X_{t})\mathbb{\ E}\
(X_{u}X_{v})(T-u)^{2H-2}|v-t|^{2H-2}dudvdt  \\
&&\le   \frac{1}{T}%
\int_{[0,T]^{3}}(T-t)^{2H-2}|u-v|^{2H-2}(T-u)^{2H-2}|v-t|^{2H-2}dudvdt \\
&&\leq
 C_{\theta, H}^2T^{8H-6}%
\int_{[0,1]^{3}}(1-t)^{2H-2}|u-v|^{2H-2}(1-u)^{2H-2}|v-t|^{2H-2}dudvdt,
\end{eqnarray*}%
which converges to $0$ as $T$ tends to infinity when $H<\frac 34  $. Hence,  by l'H\^opital
rule, $\mathbb{E}\left( |C_T^{(1)} - \mathbb{E}(C_T^{(1)}) |^2 \right)$  converges to $0$ as $T$ tends to infinity. In the same way we can prove that  $\mathbb{E}\left( |C_T^{(i)} - \mathbb{E}(C_T^{(i)}) |^2 \right)$ converges to zero as $T$ tends to infinity, for $i=2,3$  when $H<3/4$.
By triangular inequality, we see that
\begin{eqnarray*}
&&\mathbb{  E} \left[ \left( \| DF_T \|^2_{\mathcal{H}}-\mathbb{ E} ( \| DF_T \|^2_{\mathcal{H}})\right)^2 \right]\\
&=& \mathbb{ E}  \left(  |C_T^{(1)}+C_T^{(2)}+C_T^{(3)}-\mathbb{  E}  (C_T^{(1)}+C_T^{(2)}+C_T^{(3)})|^2 \right)  \\
&\le& 9  \sum_{i=1}^3 E\left( |C_T^{(i)} - \mathbb{E}(C_T^{(i)}) |^2 \right)  \\
&\rightarrow & 0\,.
\end{eqnarray*}
Taking into account that
\[
\lim_{T\rightarrow \infty} \mathbb{ E} ( \| DF_T
\|^2_{\mathcal{H}})=2\lim_{T\rightarrow \infty} \mathbb{  E}  (F_T^2),
\]
we conclude the proof of (ii).
This completes the proof of the theorem.  \hfill $\Box$

\medskip
If one replaces the It\^o type integral in \hbox{(\ref{eq1})} by the
path-wise Riemann-Stieltjes integral, then we can obtain the following estimator
\[
\widehat \theta_T^{\prime}=-\frac{\int_0^T X_t \circ d X_t}{\int_0^T X_t^2 dt}
=  \frac{ X_T^2}{2\int_0^T X_t^2 dt},
\]
which converges to zero in $L^2$ as $T$ tends to infinity from  Lemma \ref{L1} and
(\ref{f5}).

\setcounter{equation}{0}
\section{An alternative estimator}

Suppose in this section that $H>\frac 12$. We introduce the following estimator
 \begin{equation}  \label{g11}
\widetilde \th_T:=\left(\frac{1}{\si^2H\Ga(2H) T}\int_0^T X_t^2
dt\right)^{-\frac{1}{2H}}.
\end{equation}
From \eref{e.2.1}, we see that $\tilde \th_T$ 
 converges to $\th$ almost surely as $T\rightarrow \infty$.
Theorem \ref{t.3.3} allows us to derive  the rate of
convergence in the approximation of $\theta$ by  $\widetilde \th_T$.

\begin{theorem}
Suppose $H\in \left( \frac 12, \frac 34 \right)$.  Then
\begin{equation}
\sqrt{T}\left[ \widetilde {\theta }_{T}-\theta \right] \overset{\mathcal{L}}{%
\rightarrow }  N\left(0,  \frac\th{(2H)^2} \si_H^2\right) \ ,
\end{equation}%
as $T$ tends to infinity, where $\si_H$ is defined in (\ref{si}).
\end{theorem}

\textit{Proof} \quad  
From Equation \eref{e.3.7}, we have
\[
\int_0^T X_t^2 dt=\frac{\si^2 \al_H \int_0^T \int_0^t
\xi^{2H-2}e^{-\th \xi} d\xi dt-X_T^2/2}{\widehat \th_T}\,.
\]
Thus
\begin{eqnarray*}
\sqrt T\left[ \widetilde \th_T-\th\right]
&=& \sqrt T\left[ \left(\frac{ H\Ga(2H)    }{\al_H \frac{1}{T} \int_0^T \int_0^t
\xi^{2H-2}e^{-\th \xi} d\xi dt-\frac{X_T^2}{2T} }\right)^{\frac1{2H}} \widehat \th_T ^{\frac1{2H}}-\th \right]
\end{eqnarray*}
From Lemma  \ref{L6} it follows that 
\[
\al_H \frac{1}{T} \int_0^T \int_0^t
\xi^{2H-2}e^{-\th \xi} d\xi dt-\frac{X_T^2}{2T}   =\al_H\Ga(2H-1) \th^{1-2H}+o(\frac 1{\sqrt{T}}),
\]
where $o(\frac 1{\sqrt{T}})$  denotes a random  variable  $H_T$ such that $\sqrt{T}H_T$ converges to zero  almost surely as $T$ tends to infinity.
Therefore,
\begin{eqnarray}
 \left(\frac{ H\Ga(2H)    }{\al_H \frac{1}{T} \int_0^T \int_0^t
\xi^{2H-2}e^{-\th \xi} d\xi dt-\frac{X_T^2}{2T} }\right)^{\frac1{2H}}  \notag
&=&\left(\frac{1  }{ \th^{1-2H}+o(\frac1{\sqrt{T}}) }\right)^{\frac1{2H}}\\
&=& \th^{1-\frac1{2H}}+o(\frac1{\sqrt{T}}) \,.  \label{g2}
\end{eqnarray}
On the other hand, we can write
\[
\sqrt T\left[   \widehat \th_T ^{\frac1{2H}}-\th^{\frac1{2H}} \right]
 =  \sqrt T \left[ \frac1{2H}\th^{\frac1{2H}-1} \left(\widehat \th_T-\th\right)+ \frac 12 \left(  \widehat \th_T-\th \right)^2 \theta^*_T\right],
\]
where  $\theta^*_T$ is a random point between $\theta$ and  $\widehat \th_T$. From Theorem \ref{t.3.3} we obtain the following convergence in law as $T$ tends to infinity:
\begin{equation}
\sqrt T\left[   \widehat \th_T ^{\frac1{2H}}-\th^{\frac1{2H}} \right]  \label{g3}
 \rightarrow   N\left(0, \frac1{(2H)^2}\th^{\frac1{2H}} \si_H^2\right)\,.
\end{equation}
Finally, from the decomposition
\begin{eqnarray*}
\sqrt T\left[ \widetilde \th_T-\th\right]
&=& \sqrt T\left[ \left(\frac{ H\Ga(2H)    }{\al_H \frac{1}{T} \int_0^T \int_0^t
\xi^{2H-2}e^{-\th \xi} d\xi dt-\frac{X_T^2}{2T} }\right)^{\frac1{2H}} -\th^{1-\frac1{2H}}
\right] \widehat \th_T ^{\frac1{2H}}\\
&& + \sqrt T  \th^{1-\frac1{2H}} \left[
 \widehat \th_T ^{\frac1{2H}} -\th ^{\frac1{2H}}\right], 
\end{eqnarray*}
and using (\ref{g2}) and (\ref{g3}) we deduce the desired convergence.
\hfill $\Box$

 \setcounter{equation}{0}
\section{Appendix}

In the sequel we  present some calculations used in the paper.
\begin{lemma} \label{A1}
For any   $H\in \left( \frac 12, 1 \right)$
\[
(2H-1) \int_0^\infty \int_0^\infty e^{-(s+u)}   |u-s|^{2H-2} duds=  \Gamma(2H).
 \]
 \end{lemma}

\medskip
\noindent
\textit{Proof} \quad
We can write, by the change-of-variables $u-s=x$,
\begin{eqnarray*}
&& \int_0^\infty \int_0^\infty e^{-(s+u)}   |u-s|^{2H-2} dsdu
 =2 \int_0^\infty \int_0^u e^{-(s+u)}   (u-s)^{2H-2}  dsdu   \\
 &&=2 \int_0^\infty \int_0^u e^{-2u+x}   x^{2H-2}  dxdu.
 \end{eqnarray*}
Integrating first in the variable $u$  and using that $(2H-1)\Gamma(2H-1) =\Gamma(2H)$ we conclude the proof.
\hfill $\Box$

\begin{lemma} \label{L6}
Let   $Y_t$ be  the stationary Gaussian process defined in (\ref{f6}), where $H>1/2$.  
Then, for any $\alpha>0$,  $\frac{Y_T}{T^\alpha}$ converges almost surely to zero as $T$ tends to infinity.
\end{lemma}

\medskip
\noindent
\textit{Proof} \quad
The covariance of the process is $Y_t$ is, using Lemma \ref{A1} to compute $\mathrm{Var}(\xi)$,
\begin{eqnarray*}
 \mathrm{Cov}(Y_0, Y_t)
&=& e^{-\th t}\mathbb{E}\left(\xi \left[\xi +\si \int_0^t e^{ \th u} dB_u^H\right]\right)\\
&=& e^{-\th t} \left[\mathrm{Var}(\xi ) +\si^2   \alpha_H   \int_0^t \int_{-\infty} ^0 e^{\th u+\th v  }|u-v|^{2H-2} du dv
   \right] \\
&=& e^{-\th t} \left[\mathrm{Var}(\xi ) +\si^2   \alpha_H  \int_0^t \int_{v} ^\infty   e^{-\th x+2 \th v  } x^{2H-2} dx dv
   \right] \\
&=& e^{-\th t} \left[\si^2\th^{-2H}H\Gamma(2H) +\si^2     \{\ \th^{1-2H} H\Gamma(2H ) t-\frac12  t^{2H}+ o(t^{2H})\ \}
   \right] \\
   &=& \si^2\th^{-2H}H \Gamma(2H)  \left[ 1   - \frac{\th^{2H} }
{ \si^2 \Gamma(2H+1) } t^{2H}+ o(t^{2H})\ \}
   \right]\,.  \\
\end{eqnarray*}
 Then the result lemma from Theorem 3.1 of Pickands  \cite{Pi}.
 \hfill $\Box$

\begin{lemma}
\label{A2}  Let $I_T$ given by (\ref{e5}). When $\frac 12<H<\frac 34$, we have
\begin{equation}
\lim_{T\rightarrow \infty }\frac{I_{T}}{T}=    \theta  %
 ^{1-4H}  \gamma_H,
\label{e.3.5}
\end{equation}
where
\[
\gamma_H=\left( 8H-2\right) \Gamma (2H-1)^{2}
+(16H-4) \frac{\Gamma(2H-1)\Gamma(3-4H)\Gamma(4H-2)}{\Gamma(2-2H)}.
\]
\end{lemma}

\medskip
\noindent
\textit{Proof} \quad
Taking the derivative with respect to $T$ we have
\begin{equation}
\frac{dI_{T}}{dT}= 4\int_{[0,T]^{3}}e^{{-\theta }\left( T-u_{2}\right)
-{\theta }|s_{1}-u_{1}|}\left( T-s_{1}\right)
^{2H-2}|u_{2}-u_{1}|^{2H-2}du_{1}du_{2}ds_{1}.  \label{e.3.6}
\end{equation}%
Making the change of variable $T-u_{2}=x_{1}$, $T-u_{1}=x_{2}$, and $%
T-s_{1}=x_{3}$ yields
\begin{equation*}
\frac{dI_{T}}{dT}=4 \int_{[0,T]^{3}}e^{{-\theta }x_{1}-{\theta }%
|x_{2}-x_{3}|}x_{3}^{2H-2}|x_{1}-x_{2}|^{2H-2}dx_{1}dx_{2}dx_{3}.
\end{equation*}%
As a consequence,%
\begin{equation*}
\lim_{T\rightarrow \infty }\frac{dI_{T}}{dT}=4 \int_{[0,\infty
)^{3}}e^{{-\theta }x_{1}-{\theta }%
|x_{2}-x_{3}|}x_{3}^{2H-2}|x_{1}-x_{2}|^{2H-2}dx_{1}dx_{2}dx_{3},
\end{equation*}%
and this integral is finite. Indeed, we can decompose this integral into the
integrals in the six disjoint regions $\{x_{\sigma (1)}<x_{\sigma
(2)}<x_{\sigma (3)}\}$, where $\sigma $ runs over all permutations of the
indices $\{1,2,3\}$. In the case $x_{1}<x_{3}<x_{2}$ making the change of
variables $x_{1}=a$, $x_{3}-x_{1}=b$, and $x_{2}-x_{3}=c$, we obtain%
\begin{eqnarray*}
&&\int_{\lbrack 0,\infty )^{3}}e^{-\theta (a+c)}(a+b)^{2H-2}(b+c)^{2H-2}dadbdc
\\
&\leq &\int_{[0,\infty )^{3}}e^{-\theta (a+c)}b^{4H-4}dadbdc,
\end{eqnarray*}%
which is finite because $H<\frac{3}{4}$. The other cases are simpler and can
be handled in a similar way.
We can write
\begin{equation*}
\int_{\lbrack 0,\infty )^{3}}e^{{-\theta }x_{1}-{\theta }%
|x_{2}-x_{3}|}x_{3}^{2H-2}|x_{1}-x_{2}|^{2H-2}dx_{1}dx_{2}dx_{3}=
 \theta   ^{1-4H}d_H,
\end{equation*}%
 where
\begin{equation} \label{eq5}
d_H=\int_{[0,\infty )^{3}}e^{-x-|y-z|}z^{2H-2}|x-y|^{2H-2}dxdydz.
\end{equation}
The integral  in (\ref{eq5})  can be simplified as follows. First we make
the change of variables $y\mapsto w$, where $w=y-x$, and we obtain
\begin{eqnarray*}
d_H &=&\int_{0}^{\infty }\int_{0}^{\infty }\int_{-x}^{\infty
}e^{-x-|x+w-z|}z^{2H-2}|w|^{2H-2}dwdxdz \\
&=&\int_{0}^{\infty }\int_{0}^{\infty }\int_{ z-x }^{\infty }\
e^{-(2x+w-z)}z^{2H-2}|w|^{2H-2}dwdxdz \\
&&+\int_{0}^{\infty }\int_{0}^{\infty }\int_{-x\ }^{z-x}\
e^{-(z-w)}z^{2H-2}|w|^{2H-2}dwdxdz.
\end{eqnarray*}%
Integrating in $x$ we get%
\begin{eqnarray*}
d_H &=&\frac{1}{2}\int_{0}^{\infty }\int_{-\infty }^{\infty }\ e^{-2\left[
  (z-w)\vee 0\right] -(w-z)}z^{2H-2}|w|^{2H-2}dwdz \\
&&+\int_{0}^{\infty }\int_{-\infty }^{\infty }\left[ (z-w)-\left( (-w)\vee
0\right) \right] _{+}\ e^{-(z-w)}z^{2H-2}|w|^{2H-2}dwdz.
\end{eqnarray*}%
Therefore,%
\begin{eqnarray*}
d_H &=&\frac{1}{2}\int_{0}^{\infty }\int_{0}^{\infty }\ e^{-2\left[ (z-w)\vee 0%
\right] -(w-z)}z^{2H-2}w^{2H-2}dwdz \\
&&+\frac{1}{2}\int_{0}^{\infty }\int_{0}^{\infty }\
e^{-(z+w)}z^{2H-2}w^{2H-2}dwdz \\
&&+\int_{0}^{\infty }\int_{0}^{\infty }\left[ (z-w)\right] _{+}\
e^{-(z-w)}z^{2H-2}w^{2H-2}dwdz \\
&&+\int_{0}^{\infty }\int_{0}^{\infty }\ e^{-(z+w)}z^{2H-1}w^{2H-2}dwdz,
\end{eqnarray*}%
and  we obtain%
\begin{equation}  \label{f1}
d_H  = f_H + \left( 2H-\frac{1}{2}\right) \Gamma (2H-1)^{2},
\end{equation}
where
\[
f_H= \int_{0}^{\infty } \int_{0}^ z (1+z-w)    e^{-(z-w)}  z^{2H-2}w^{2H-2}dwdz.
\]
Making the change-of-variables $z-w=x$ yields
\[
f_H=  \int_{0}^{\infty } \int_{0}^{\infty }  (1+x)    e^{-x}  (w+x)^{2H-2}w^{2H-2}dwdx.
\]
Substituting the equality $  (w+x)^{2H-2}=\frac 1{\Gamma(2-2H)}  \int_{0}^{\infty }
e^{-\xi (w+x)} \xi^{1-2H} d\xi$ in $f_H$ we obtain
\begin{eqnarray}
f_H&=& \frac 1{\Gamma(2-2H)} \int_{0}^{\infty } \int_{0}^{\infty } \int_{0}^{\infty } (1+x)    e^{-x-\xi(w+x)}   w^{2H-2}  \xi^{1-2H} d\xi dwdx  \notag\\
&&=  \frac{\Gamma(2H-1)}{\Gamma(2-2H)} \int_{0}^{\infty }  \int_{0}^{\infty } (1+x)    e^{-x-\xi x}      \xi^{2-4H} d\xi  dx \notag\\
&&= \frac{\Gamma(2H-1) \Gamma(3-4H)}{\Gamma(2-2H)} \int_{0}^{\infty }    (1+x)    e^{-x  }        x^{4H-3}  dx  \notag\\
&&=(4H-1) \frac{\Gamma(2H-1)\Gamma(3-4H)\Gamma(4H-2)}{\Gamma(2-2H)}.
\label{f2}
\end{eqnarray}
Finally from (\ref{f1}) and (\ref{f2}) we get the desired result.
\hfill $\Box$

 \begin{lemma}  \label{lem1}
Let $X_{t}$ be given by (\ref{eq1}). We have%
\begin{equation}
\mathbb {E}  \left[ \int_{s}^{T}e^{-\theta (\xi -s)}dB_{\xi
}^{H}\int_{t}^{T}e^{-\theta (\eta -t)}dB_{\eta }^{H}\right] \leq
C_{\theta,H}|t-s|^{2H-2},  \label{e.3.22}
\end{equation}%
and%
\begin{equation}
\mathbb{E}  \left[ X_{t}X_{s}\right] \leq  \sigma^2C_{\theta,H}|t-s|^{2H-2},  \label{e.3.23}
\end{equation}%
for some constant $C_{\theta, H}>0$.
\end{lemma}

\medskip
\noindent
\textit{Proof} \quad
Let us assume that $s<t$.  We can write using (\ref{eq9})
\begin{eqnarray*}
&&\mathbb{ E}  \left[ \int_{s}^{T}e^{-\theta (\xi -s)}dB_{\xi
}^{H}\int_{t}^{T}e^{-\theta (\eta -t)}dB_{\eta }^{H}\right]  \\
& &\quad =  \alpha_H\int_{t}^{T}\int_{s}^{T}e^{-\theta (\xi -s)}e^{-\theta (\eta -t)}|\xi -\eta
|^{2H-2}d\xi d\eta  = \alpha_H (B_T^{(1)} + B_T^{(2)}),
\end{eqnarray*}
where
\[
B_T^{(1)}=\int_{t}^{T}\int_{t}^{T}e^{-\theta (\xi -s)}e^{-\theta (\eta -t)}|\xi -\eta
|^{2H-2}d\xi d\eta
\]
and
\[
B_T^{(2)}=\int_{t}^{T}\int_{s}^{t}e^{-\theta (\xi -s)}e^{-\theta
(\eta -t)}|\xi -\eta |^{2H-2}d\xi d\eta .
\]
It is easy to see that $B_T^{(1)}$  is bounded by $C_{\theta,H}e^{{-\theta }|t-s|}$.
The second term can be estimated as follows
\begin{eqnarray*}
 B_T^{(2)}
&=&\int_{s}^{t}\int_{t-\xi }^{T-\xi }e^{-\theta (\xi -s+y+\xi
-t)}y^{2H-2}dyd\xi  \\
&=&\int_{0}^{T-s}y^{2H-2}dy\int_{(t-y)\vee s}^{(T-y)\wedge t}e^{-\theta
(y+2\xi -s-t)}d\xi  \\
&\leq &\frac{1}{2{\theta }}\int_{0}^{T-s}e^{{-\theta }(y+2\,(t-y)\vee
s-s-t)}y^{2H-2}dy \\
&=&\frac{1}{2{\theta }}\int_{t-s}^{T-s}e^{{-\theta }(y+s-t)}y^{2H-2}dy+\frac{%
1}{2{\theta }}\int_{0}^{t-s}e^{{-\theta }(y+s-t)}y^{2H-2}dy \\
&\leq &C_{\theta,H}|t-s|^{2H-2}\int_{t-s}^{T-s}e^{{-\theta }(y+s-t)}dy+C_{\th %
}\int_{0}^{t-s}y^{2H-2}dy \\
&\leq &C_{\theta,H}|t-s|^{2H-2}\,.
\end{eqnarray*}%
This proves (\ref{e.3.22}). The inequality (\ref{e.3.23}) can be proved in a
similar way (see also \cite{CKM}).
\hfill $\Box$

\end{document}